\documentclass{article}

\usepackage{amsmath} 
\usepackage{amssymb}
\usepackage{amsthm}
\usepackage{enumerate}

\usepackage{tikz}

\theoremstyle{definition}
\newtheorem{defi}{Def.}

\theoremstyle{plain}
\newtheorem{lemma}[defi]{Lemma}
\newtheorem{theorem}[defi]{Theorem}

\newtheorem{kor}[defi]{Corollary}

\theoremstyle{remark}

\newcommand{\R}[1]{\mathbb{R}^{#1}}
\newcommand{\prf}[1][]{\noindent {\bf Proof{#1}:} }

\DeclareMathOperator{\bd}{bd}

\DeclareMathOperator{\Cov}{Cov}

\DeclareMathOperator{\dist}{dist}
\DeclareMathOperator{\Var}{Var}
\DeclareMathOperator{\Lip}{Lip}

\title{A central limit theorem for Lebesgue integrals of random fields}
\author{J\"urgen Kampf}

\begin{document}
\maketitle

\begin{abstract}

In this paper we show a central limit theorem for Lebesgue integrals of stationary $BL(\theta)$-dependent random fields as the integration domain grows in Van Hove-sense. Our method is to use the (known) analogue result for discrete sums. As applications we obtain various multivariate versions of this central limit theorem. 

\end{abstract}

\section{Introduction}
Random fields are collections of random variables indexed by the Euclidean space $\mathbb{R}^d$. They have applications in various branches of science, e.g.\ in medicine \cite{ATW, TW}, in geostatistics \cite{ChDe, Wack} or in material science \cite{McSt, Tor}.  

The aim of the present paper is to establish a central limit theorem for integrals $\int_{W_n} X(t)\, dt$, where $(W_n)_{n\in\mathbb{N}}$ is a sequence of compact subsets of $\R{d}$ and $(X(t))_{t\in\mathbb{R}^d}$ is a random field. The sequence $(W_n)_{n\in\mathbb{N}}$ of integration domains is assumed to \emph{grow in Van Hove-sense (VH-sense)}, i.e.\ 
\[ \lim_{n\to\infty} \lambda_d((\bd W_n)+B^d) / \lambda_d(W_n) = 0, \] 
where $\lambda_d$ denotes the Lebesgue measure, $\bd W$ is the boundary of $W\subseteq \R{d}$, $A+B:=\{ a+b \mid a\in A,\ b\in B\}$ for two subsets $A,B\subseteq \R{d}$ and $B^d:=\{x\in\R{d} \mid \|x\|\le 1 \}$ is the closed Euclidean unit ball. 

%We will now define the assumption of $BL(\theta)$-dependence which the random field  $(X(t))_{t\in\mathbb{R}^d}$ has to fulfill.  

%\begin{defi}
%Let $\theta=(\theta_r)_{r\in \mathbb{N}}$ be a monotonically decreasing sequence with $\lim_{r \to \infty} \theta_r=0$.
% An $\mathbb{R}^s$-valued random field $(X(t))_{t\in \mathbb{R}^{d}}$ is called \emph{$BL(\theta)$-dependent} if for any $\Delta>1$ and any disjoint, finite sets $I,J\subseteq T(\Delta)$ with $\dist(I,J)\ge r$ and all bounded Lipschitz functions $f:\mathbb{R}^{s\cdot \# I}\to \mathbb{R}$ and $g:\mathbb{R}^{s\cdot \# J}\to \mathbb{R}$ we have 
%\[ \Cov(f(X_I), g(X_J)) \le \Psi(\# I, \# J, f, g)\Delta^d \theta_r. \] 
%\end{defi}

%For example, all associated random fields whose covariance function is decreasing fast enough are $BL(\theta)$-dependent. 

The main result of the present paper is the following (the notion of $BL(\theta)$-dependence will be defined in Subsection \ref{ss:ass}). 

\begin{theorem}\label{T:uCLT}
Let $\theta=(\theta_r)_{r\in\mathbb{N}}$ be a monotonically decreasing zero sequence. Let $(X(t))_{t\in \R{d}}$ be a measurable, stationary, $BL(\theta)$-dependent $\mathbb{R}$-valued random field such that 
\[ \int_{\R{d}}| \Cov\big( X(0), X(t) \big)|\, dt < \infty.\]
Let $(W_n)_{n\in \mathbb{N}}$ be a VH-growing sequence of subsets of $\R{d}$. 
Then
\[\frac{ \int_{W_n} X(t)\, dt  - \mathbb{E}\, X(0) \lambda_d(W_n)}{\sqrt{\lambda_d(W_n)}}\to \mathcal{N}(0,\sigma^2), \quad n\to \infty,\]
in distribution, where
\[ \sigma^2 := \int_{\R{d}} \Cov\big( X(0), X(t) \big)\, dt.\] 
\end{theorem}

There is a wide literature on similar results, where mixing conditions are assumed instead of $BL(\theta)$-dependence, see e.g.\ \cite{BulZ, Gor, IL, Leon}. For $BL(\theta)$-dependent random fields there are no central limit theorems for Lebesgue integrals up to now. However, there are such results for discrete sums \cite{BS} and for Lebesgue measures of excursion sets \cite{BST}. In the latter paper the random field is in fact assumed to be quasi-associated, which is a slightly stronger assumption than $BL(\theta)$- dependence.   

This paper is organized as follows: In Section \ref{s:Prep} we collect preliminaries about associated random variables, random fields and functions of bounded variation. Section \ref{s:Uni} is devoted to the proof of the main theorem. In Section \ref{s:Multi} we present several examples how the main result can be extended to a multivariate central limit theorem. The case that the random field is of the form $(f(X(t)))_{t\in\mathbb{R}^d}$ for some deterministic function $f:\mathbb{R} \to\mathbb{R}^s$ and some random $\mathbb{R}$-valued field $(X(t))_{t\in\mathbb{R}^d}$ will be of particular interest.

\section{Preliminaries}\label{s:Prep}
 
\subsection{Association concepts}\label{ss:ass}
In this subsection we introduce different association concepts and discuss their relations. 

We start with the broadest appearing in this paper, namely $BL(\theta)$-dependence.

For finite subsets $I, J\subseteq \mathbb{R}^d$ we put $\dist(I, J):= \min\{\|x-y\|_1: x\in I, \, y \in J\}$, where $\|\cdot\|_1$ is the $\ell_1$-norm. For two Lipschitz functions $f:\mathbb{R}^{n_1}\to \mathbb{R}$ and $g:\mathbb{R}^{n_2}\to \mathbb{R}$ we put
\[ \Psi(n_1, n_2, f, g) = \min\{n_1, n_2\}\Lip(f)\Lip(g),\]
where
\[ \Lip(f):=\sup\Big\{\frac{|f(x)-f(y)|}{\|x-y\|_1} \mid x, y \in \mathbb{R}^n, \ x \ne y\Big\} \]
denotes the (optimal) Lipschitz constant of a Lipschitz function $f: \mathbb{R}^n \to \mathbb{R}$.

For a random field $(X(t))_{t\in \R{d}}$, a finite subset $I=\{t_1,\dots,t_n\}\subseteq \R{d}$ with $n$ elements and a function $f$ on $\R{n}$ we abbreviate $f(X_I):=f(X(t_1), \dots, X(t_n))$. If such an abbreviation $X_I$ appears more than once within one formula, then always the same enumeration of the elements of $I$ has to be used.   

For a set $M$ let $\# M$ denote the number of elements of $M$.

 Furthermore, for $\Delta>0$ we put
\[ T(\Delta):= \{ (j_1/\Delta, \dots, j_d/\Delta) \mid (j_1,\dots, j_d)\in \mathbb{Z}^d\}. \]

\begin{defi}
Let $\theta=(\theta_r)_{r\in \mathbb{N}}$ be a monotonically decreasing sequence with $\lim_{r \to \infty} \theta_r=0$.
\begin{enumerate}[(i)]
\item An $\mathbb{R}^s$-valued random field $(X(t))_{t\in \mathbb{R}^{d}}$ is called \emph{$BL(\theta)$-dependent} if for any $\Delta>1$ and any disjoint, finite sets $I,J\subseteq T(\Delta)$ with $\dist(I,J)\ge r$ and all bounded Lipschitz functions $f:\mathbb{R}^{s\cdot \# I}\to \mathbb{R}$ and $g:\mathbb{R}^{s\cdot \# J}\to \mathbb{R}$ we have
\[ \Cov(f(X_I), g(X_J)) \le \Psi(\# I, \# J, f, g)\Delta^d \theta_r. \] 
\item An $\mathbb{R}^s$-valued random field $(X(t))_{t\in \mathbb{Z}^d}$ is called \emph{$BL(\theta)$-dependent} if for any disjoint, finite sets $I,J\subseteq \mathbb{Z}^{d}$ with $\dist(I,J)\ge r$ and all bounded Lipschitz functions $f:\mathbb{R}^{s\cdot \# I}\to \mathbb{R}$ and $g:\mathbb{R}^{s\cdot \# J}\to \mathbb{R}$ we have
\[ \Cov(f(X_I), g(X_J)) \le \Psi(\# I, \# J, f, g) \theta_r. \]
\end{enumerate}  
\end{defi}

\begin{lemma}\label{L:BLlim}
Let $\theta=(\theta_r)_{r\in \mathbb{N}}$ be a monotonically decreasing sequence with $\lim_{r \to \infty} \theta_r=0$. For $T=\mathbb{Z}^d$ or $T=\mathbb{R}^d$, let $(X^{(n)}(t))_{t \in T}, n\in \mathbb{N},$ be a sequence of $BL(\theta)$-dependent random fields such that the finite-dimensional distributions converge to those of a field $(X(t))_{t\in T}$. Then $(X(t))_{t\in T}$ is also $BL(\theta)$-dependent.
\end{lemma}
\prf By the definition of convergence in probability we get $\lim_{n\to\infty} \mathbb{E}\, f(X_I^{(n)})g(X_J^{(n)})=\mathbb{E}\, f(X_I)g(X_J)$, $\lim_{n\to\infty} \mathbb{E}\, f(X_I^{(n)})=\mathbb{E}\, f(X_I)$ and $\lim_{n\to\infty} \mathbb{E}\, g(X_J^{(n)})=\mathbb{E}\, g(X_J)$ for any finite sets $I,J\subseteq T$ and bounded Lipschitz continuous functions $f:\mathbb{R}^{\#I}\to \mathbb{R}$ and $g:\mathbb{R}^{\#J}\to \mathbb{R}$, which yields the assertion. \qed

\begin{lemma}\label{L:LipBL}
Let $\theta=(\theta_r)_{r\in \mathbb{N}}$ be a monotonically decreasing sequence with $\lim_{r \to \infty} \theta_r=0$. Let $(X(t))_{t\in T}$ be a $BL(\theta)$-dependent random field and let $f:\mathbb{R}^{s} \to \mathbb{R}^{s'}$ be a Lip\-schitz function. Then there is a monotonically decreasing sequence $\theta'=(\theta'_r)_{r\in \mathbb{N}}$ with $\lim_{r \to \infty} \theta'_r=0$ such that $(f(X(t)))_{t\in T}$ is $BL(\theta')$-dependent.
\end{lemma}
\prf For a function $f: V\to W$ and a finite set $I$ let $f_I$ denote the function $V^{\# I} \to W^{\# I}$, $(x_1,\dots, x_{\# I}) \mapsto (f(x_1),\dots, f(x_{\# I}))$.  We put $\theta'_r:=\Lip(f)^2\cdot \theta_r$. Let $I,J$ be two disjoint finite sets with $\dist(I,J)\ge r$ and let $\tilde f:\mathbb{R}^{s'\cdot\# I}\to \mathbb{R}$ and $g:\mathbb{R}^{s'\cdot\# J}\to \mathbb{R}$ be two Lip\-schitz functions. Then
\begin{align*}
\Cov(\tilde f(f_I(X_I)), g(f_J(X_J)) ) &\le \min\{\#I, \#J\}\cdot \Lip(\tilde f \circ f_I)\cdot \Lip(g\circ f_J)\cdot\Delta^d \cdot\theta_r \\
&\le \min\{\#I, \#J\}\cdot \Lip(\tilde f)\cdot \Lip(g)\cdot\Delta^d \cdot\theta'_r.\qquad\qed 
\end{align*}
%where $\theta'_r:=\theta_r\cdot Lip(f)^2$.\qed

%\begin{lemma}\label{L:BLlim}
%Let $(X^{(n)}(t))_{t \in T}, n\in \mathbb{N},$ be a sequence of QA random fields such that the finite-dimensional distributions converge to those of a field $(X(t))_{t\in T}$. Then $(X(t))_{t\in T}$ is also QA. \prf Let $I, J\subseteq T$ be two finite sets and $f:\mathbb{R}^{\# I}\to\mathbb{R}$ and $g:\mathbb{R}^{\# J}\to\mathbb{R}$ two Lipschitz functions. Then the Continuous Mapping Theorem implies that $f(X_I^{(n)}) \to f(X_I)$ and $g(X_J^{(n)}) \to g(X_J)$ in distribution. So the Dominated Convergence Theorem implies $\lim_{n\to \infty} \Cov(f(X_I^{(n)}), g(X_J^{(n)})) = \Cov(f(X_I), g(X_J))$, which yields the assertion. 
%\end{lemma}

An $\R{s}$-valued random field $(X(t))_{t\in T}$ is called \emph{positively associated} (PA) if
\[ \Cov(f(X_I), g(X_J)) \ge 0 \]
for any finite sets $I, J  \subseteq T$ and functions $f:\mathbb{R}^{s\cdot \# I}\to \mathbb{R}$ and $g:\mathbb{R}^{s\cdot \# J}\to\mathbb{R}$ which are bounded and monotonically increasing in every coordinate.

For a Lipschitz function $f:\mathbb{R}^n\to\mathbb{R}$ we define coordinate-wise Lipschitz constants by
\begin{align*}
\Lip_{k}(f) = \sup \Big\{\frac{|f(x_1,\dots, x_{k-1}, y_k, x_{k+1},\dots, x_n) - f(x_1,\dots, x_{k-1}, z_k, x_{k+1},\dots, x_n)|}{|y_k-z_k|} &\mid \\
x_1,\dots, x_{k-1}, x_{k+1},\dots, x_n, y_k, z_k\in\mathbb{R},\ y_k\ne z_k &\Big\}, \ \, k\in\{1,\dots, n\}.
\end{align*} 
An $\R{s}$-valued random field $(X(t))_{t\in T}$ with $\mathbb{E}[X_{k}(t)^2]<\infty, t\in T, k=1,\dots, s,$ is called \emph{quasi-associated} (QA) if
\[ |\Cov(f(X_I), g(X_J))| \le \sum_{t\in I} \sum_{k=1}^s \sum_{u\in J} \sum_{l=1}^s \Lip_{t,k}(f)\cdot \Lip_{u,l}(g) |\Cov(X_{k}(t), X_{l}(u))| \]
for any finite sets $I, J  \subseteq T$ and Lipschitz continuous functions $f:\mathbb{R}^{s\cdot \# I}\to \mathbb{R}$ and $g:\mathbb{R}^{s\cdot \# J}\to\mathbb{R}$. 

It is well known that every PA random field is also QA, see e.g.\ Theorem 5.3 in \cite[p.\ 89]{BS} (this theorem is only formulated in the special case $s=1$ and $T=\mathbb{Z}^d$, but the proof holds in the present setting). 

%The following lemma is clear since the composition of two Lipschitz continuous functions is again Lipschitz continuous with Lipschitz constant not bigger than the product of the Lipschitz constants of the two individual functions. 

%\begin{lemma}
%Let $(X(t))_{\t\in\mathbb{R}^d}$ be an $\mathbb{R}^s$-valued QA random field and    

\begin{lemma}\label{l:QABl}
Let $(X(t))_{t\in\mathbb{R}^d}$ be an $\mathbb{R}^s$-valued QA random field. Assume that there are $c>0$ and $\epsilon>0$ with
\[ \Cov( X_i(t_1), X_j(t_2) ) \le c \cdot \|t_1-t_2\|^{-d-\epsilon}_\infty\]
for $t_1,t_2\in \R{d}$ and $i,j=1,\dots, s$.
Then $ (X(t))_{t\in\mathbb{R}^d}$ is $BL(\theta)$-dependent for some monotonically decreasing zero sequence $\theta$. 
\end{lemma}
\prf Let $r>0$, $\Delta>1$ and let $I,J\subseteq T(\Delta)$ be finite with $\dist(I,J)\ge r$, w.l.o.g.\ $\# I \le \# J$. Moreover, let $f:\R{s\cdot \# I}\to\R{}$ and $g:\R{s\cdot \# J} \to \R{}$ be bounded and Lipschitz continuous. Then

\begin{align*}
\Cov(f(X_I), g(X_J)) &\le \sum_{t\in I} \sum_{k=1}^s \sum_{u\in J} \sum_{l=1}^s \Lip_{t,k}(f)\cdot \Lip_{u,l}(g) \Cov(X_{k}(t), X_{l}(u)) \\
&\le s^2\cdot \# I \cdot \Lip(f)\cdot \Lip(g) \cdot \max_{t,k,l} \sum_{u\in J} \Cov(X_{k}(t), X_{l}(u))\\
&\le s^2\cdot \# I \cdot \Lip(f)\cdot \Lip(g) \cdot \max_{t\in I} \sum_{u\in J} c \cdot \|t-u\|^{-d-\epsilon}_\infty.
\end{align*}
We have for fixed $t\in I$, if $r>1$,
\begin{align*} 
\sum_{u\in J} \|t-u\|^{-d-\epsilon}_\infty %&= \sum_{j\in T(\Delta)\setminus \{0\} } \|j\|^{-d-\epsilon}_\infty \\
&\le \sum_{s=\lceil r \Delta \rceil}^\infty \big(\frac{s}{\Delta}\big)^{-d-\epsilon}\cdot \#\{v\in T(\Delta) \mid \|v\|_\infty=\frac{s}{\Delta}\} \\
&=\sum_{s=\lceil r \Delta \rceil}^\infty \big(\frac{s}{\Delta}\big)^{-d-\epsilon}\cdot \big( (2s+1)^d-(2s-1)^d \big)\\
&= \Delta^{d+\epsilon} \sum_{s=\lceil r \Delta \rceil}^\infty \sum_{\iota=0}^{d-1} s^{-d-\epsilon} {d \choose \iota} (1+(-1)^{d-\iota-1}) (2s)^\iota\\
& \le \Delta^{d+\epsilon} \int_{\lceil r \Delta \rceil-1}^\infty \sum_{\iota=0}^{d-1} {d \choose \iota} (1+(-1)^{d-\iota-1})2^\iota s^{-d-\epsilon+\iota}\, ds\\
& \le \Delta^{d+\epsilon} \sum_{\iota=0}^{d-1}  {d \choose \iota} (1+(-1)^{d-\iota-1})2^\iota \frac{(r\Delta-\Delta)^{-d-\epsilon+\iota+1}}{d+\epsilon-\iota-1}\\
& \le \Delta^{d} \sum_{\iota=0}^{d-1}  {d \choose \iota} (1+(-1)^{d-\iota-1})2^\iota \frac{(r-1)^{-d-\epsilon+\iota+1}}{d+\epsilon-\iota-1}.
\end{align*}
Putting
\[ \theta_r:= \begin{cases} c \cdot s^2\sum_{\iota=0}^{d-1}  {d \choose \iota} (1+(-1)^{d-\iota-1})2^\iota \frac{(r-1)^{-d-\epsilon+\iota+1}}{d+\epsilon-\iota-1} & \mbox{for } r>1,\\
															\frac{3^d}{\Delta^d}\cdot\max_{i=1,\dots, s} \Var(X_i(0)) +\theta_2 & \mbox{for } r=1,
							\end{cases} \]
we obtain
\[ \Cov(f(X_I), g(X_J)) \le \min\{ \# I, \# J \} \cdot \Lip(f)\cdot \Lip(g) \cdot \Delta ^d \theta_r. \qed  \]

\subsection{Random fields}

After having introduced the association concepts in subsection \ref{ss:ass}, we will now collect various other preliminaries concerning random fields.

%A sequence $(W_n)_{n\in\mathbb{N}}$ of compact subsets of $\R{d}$ is called \emph{Van Hove-growing} (VH-growing) if
%\[ \lim_{n\to\infty} \lambda_d((\bd W_n)+B^d)/\lambda_d(W_n) =0, \]
%where $\bd W$ denotes the boundary of $W$. 

The following theorem (see \cite[Ch.\ III, \S\ 3]{GS} and \cite[Prop.\ 3.1]{Roy}) says that for stationary random fields stochastic continuity and measurability are essentially equivalent.
\begin{theorem}\label{T:cont=meas}
\begin{enumerate}[(i)]
\item Let $(X(t))_{t\in\R{d}}$ be a stochastically continuous random field. Then there is a measurable modification of $(X(t))_{t\in\R{d}}$.
\item Let $(X(t))_{t\in\R{d}}$ be a stationary and measurable random field. Then $(X(t))_{t\in\R{d}}$ is stochastically continuous.
\end{enumerate}
\end{theorem}

\begin{lemma}\label{L:contin}
Let $(X_t)_{t\in\mathbb{R}^d}$ be a stationary, stochastically continuous random field with $\mathbb{E}\, X(0)^j<\infty$ for $j>0$. Then $(X_t)_{t\in\mathbb{R}^d}$ is continuous in $j$-mean. 
\end{lemma}
\prf Let $(t_n)_{n\in\mathbb{N}}$ be a sequence of points in $\R{d}$ converging to a point $t\in\mathbb{R}^d$. Then
\[\lim_{n\to\infty} \mathbb{E}\, |X(t_n)-X(t)|^j =  \lim_{n\to\infty} \int_0^\infty \mathbb{P}(|X(t_n) - X(t)|^j \ge x) \, dx = \int_0^\infty \lim_{n\to\infty} \mathbb{P}(|X(t_n) - X(t)| \ge \sqrt[j]{x}) \, dx =0.\]
We have been allowed to interchange limit and integral, since
\[\mathbb{P}(|X(t_n) - X(t)| \ge \sqrt[j]{x}) \le \mathbb{P}(|X(t_n)| \ge \tfrac{\sqrt[j]{x}}{2}) + \mathbb{P}(|X(t)| \ge \tfrac{\sqrt[j]{x}}{2}) = 2 \mathbb{P}(|X(t)| \ge \tfrac{\sqrt[j]{x}}{2})\]
due to the stationarity and 
\[\int_0^\infty 2 \mathbb{P}(|X(t)| \ge \tfrac{\sqrt[j]{x}}{2}) \, dx = \int_0^\infty 2 \mathbb{P}(|2\cdot X(t)|^j \ge x) \, dx = 2^{j+1} \mathbb{E}\, |X(0)|^j<\infty.\qquad\qed\] 

\begin{lemma}\label{L:sameint}
Let $(X(t))_{t\in\mathbb{R}^d}$ and $(Y(t))_{t\in\mathbb{R}^d}$ be two stochastically continuous and measurable random fields having the same distribution. Let $A_1,\dots, A_m\subseteq\mathbb{R}^d$ be bounded Borel sets. Assume that $\int_{A_i} X(t)\, dt$ is defined a.s.\ for $i=1,\dots,m$, i.e.\ not both the positive part and the negative part of these integrals are infinite. Then $\int_{A_1} Y(t)\, dt, \dots, \int_{A_m} Y(t)\,dt$ are defined a.s.\ as well and
\[ \bigg(\int_{A_1} X(t)\, dt, \dots, \int_{A_m} X(t)\,dt \bigg) \stackrel{d}{=} \bigg(\int_{A_1} Y(t)\, dt, \dots, \int_{A_m} Y(t)\,dt \bigg). \]
\end{lemma}

\prf By the Monotone Convergence Theorem, we may assume w.l.o.g.\ that there is some $N\in\mathbb{N}$ such that $X(t)\in [-N,N]$ and $Y(t)\in [-N,N]$ for all $t\in\R{d}$.

We define processes $(X^n(t))_{t\in\R{d}}$ and $(Y^n(t))_{t\in\R{d}}$ by putting
\[ X^n(t_1,\dots, t_d) := X\Big(\frac{z_1}{n}, \dots, \frac{z_d}{n}\Big), \quad \mbox{for all } t_1\in \Big[\frac{z_1}{n}, \frac{z_1+1}{n}\Big), \dots, t_d\in \Big[\frac{z_d}{n}, \frac{z_d+1}{n}\Big), \, z_1,\dots, z_d\in\mathbb{Z}.\]
We get
\begin{align} 
\Big\{ \int_{A_i} X^n(t)\, dt & \Big| i=1,\dots, m \Big\}\\
 &= \Big\{ \sum_{z_1, \dots, z_d\in\mathbb{Z}} \lambda_d\Big( A_i\cap \big[\frac{z_1}{n}, \frac{z_1+1}{n}\big)\times \dots \times \big[\frac{z_d}{n}, \frac{z_d+1}{n}\big)\Big) X\big(\frac{z_1}{n}, \dots, \frac{z_d}{n}\big) \Big| i=1,\dots, m \Big\} \notag\\
& \stackrel{d}{=} \Big\{ \sum_{z_1, \dots, z_d\in\mathbb{Z}} \lambda_d\Big( A_i\cap \big[\frac{z_1}{n}, \frac{z_1+1}{n}\big)\times \dots \times \big[\frac{z_d}{n}, \frac{z_d+1}{n}\big)\Big) Y\big(\frac{z_1}{n}, \dots, \frac{z_d}{n}\big) \Big| i=1,\dots, m \Big\} \notag\\
&= \Big\{ \int_{A_i} Y^n(t)\, dt \Big| i=1,\dots, m \Big\}. \label{e:Riem} \end{align}

For $\epsilon, \delta >0$  we get 
\begin{align*}
\mathbb{P}\Big(\sum_{i=1}^m \big|\int_{A_i} X^n(t) \, dt - \int_{A_i} X(t) \, dt \big|>\epsilon\Big) & \le \mathbb{P}\Big(\sum_{i=1}^m \int_{A_i} |X^n(t) -X(t)| \, dt >\epsilon\Big) \\
& \le \frac{\mathbb{E}\, \sum_{i=1}^m \int_{A_i} |X^n(t) -X(t)| \, dt }{\epsilon}\\
& = \frac{\sum_{i=1}^m \int_{A_i} \mathbb{E}\, |X^n(t) -X(t)| \, dt }{\epsilon}\\
& \le \frac{\sum_{i=1}^m \int_{A_i} \big(\delta + \mathbb{P}( |X^n(t) -X(t)|>\delta)\cdot 2N \big)\, dt }{\epsilon}\\
& \stackrel{n\to\infty}{\longrightarrow} \frac{\sum_{i=1}^m  \int_{A_i} \delta  \, dt }{\epsilon}\\
& = \frac{\delta \cdot \sum_{i=1}^m \lambda_d(A_i)}{\epsilon}.
\end{align*}
The limit relation holds by the Majorized Convergence Theorem, since the assumption that $(X(t))_{t\in\mathbb{R}}$ is stochastically continuous implies that $X^n(t)$ converges to $X(t)$. Since $\delta>0$ was arbitrary, we get  
\[ \mathbb{P}\Big(\sum_{i=1}^m \big|\int_{A_i} X^n(t) \, dt - \int_{A_i} X(t) \, dt \big|>\epsilon\Big) \stackrel{n\to\infty}{\longrightarrow} 0 \]
and the same way 
\[ \mathbb{P}\Big(\sum_{i=1}^m \big|\int_{A_i} Y^n(t) \, dt - \int_{A_i} Y(t) \, dt \big|>\epsilon\Big) \stackrel{n\to\infty}{\longrightarrow} 0. \]
Now (\ref{e:Riem}) yields the assertion. \qed\medskip

\subsection{Functions of bounded variation}

A function $f:\R{}\to\R{}$ is said to be of \emph{locally bounded variation} if there is a monotonically increasing function $\alpha:\mathbb{R}\to\mathbb{R}$ and a monotonically decreasing function $\beta:\mathbb{R}\to\mathbb{R}$ such that $f=\alpha+\beta$. We denote the set of such functions $\alpha$ and $\beta$ by $A$ resp.\ $B$. We put
\[f^+(x):=\begin{cases} 
\inf\{\alpha(x)\mid \alpha\in A,\, \alpha(0)=f(0)\} & \mbox{ if } x>0\\ 
 f(0) & \mbox{ if } x=0\\ 
\sup\{\alpha(x)\mid \alpha\in A,\, \alpha(0)=f(0)\} & \mbox{ if } x<0.
\end{cases}\]
It is easy to see that $f^+\in A$ and $f^-:=f-f^+\in B$. We put $h_f:=f^+-f^-$.  

\begin{lemma}\label{L:bv_decomp} Let $f:\R{} \to \R{}$ be a function of locally bounded variation.  Then $f=g\circ h_f$ for a Lipschitz continuous function $g:\R{}\to\R{}$ of Lipschitz constant $1$.
\end{lemma}
\prf For each $x\in \mathbb{R}$, for which there is $t\in\mathbb{R}$ with $h_f(t)=x$, define $g(x):=f(t)$. Now $g$ is well-defined, since for $t_1,t_2\in\mathbb{R}$ with $h_f(t_1)=h_f(t_2)$, $f$ is constant on $[t_1,t_2]$ . Clearly, $f=g\circ h_f$. Moreover, $g$ -defined on a subset of $\mathbb{R}$ so far- is Lipschitz continuous with Lipschitz constant $1$. Indeed, let $x_1,x_2\in \mathbb{R},\ x_1<x_2$, be two points for which there are $t_1,t_2\in\mathbb{R}$ with $h_f(t_1)=x_1$ and $h_f(t_2)=x_2$. Then
\[h_f(t_2)-h_f(t_1)=f^+(t_2) - f^+(t_1) - (f^-(t_2)-f^-(t_1)) \ge |f^+(t_2) - f^+(t_1) + (f^-(t_2)-f^-(t_1))| = |f(t_2) -f(t_1)|. \]
Hence $x_2-x_1\ge |g(x_2)-g(x_1)|$. 

It remains to show that $g$ has a Lipschitz continuous extension to the whole of $\mathbb{R}$. The domain of $g$ is $\mathbb{R}$ minus the union of countable many, disjoint intervals. For a point $x$ lying on the boundary of the domain of $g$ but not in the domain of $g$, choose a sequence $(x_n)_{n\in\mathbb{N}}$ such that $g(x_n)$ is defined for all $n\in\mathbb{N}$. Then $(g(x_n))_{n\in\mathbb{N}}$ is a Cauchy sequence, since $g$ is Lipschitz continuous, and hence convergent. Since $(g(x_n))_{n\in\mathbb{N}}$ is convergent for \emph{every} such sequence $(x_n)_{n\in\mathbb{N}}$, the limit is independent of the choice of the sequence. So we can put $g(x):=\lim_{n\to\infty} g(x_n)$. It is easy to see that this extension still has Lipschitz constant $1$.  Now all gaps in the domain of $g$ are \emph{open} intervals. So they can be filled by affine functions. Clearly, the Lipschitz constant is preserved again. 
\qed\medskip

\section{The univariate CLT}\label{s:Uni}

In this section we will prove Theorem \ref{T:uCLT}.

\prf For $j=(j_1,\dots,j_d)\in \mathbb{Z}^d$ we put $Q_j=\mathop{\times}_{i=1}^d [j_i,j_i+1)$ and $Z(j):=\int_{Q_j} X(t)\, dt-\mathbb{E}\, X(0)$. We will show that this random field $(Z(j))_{j\in\mathbb{Z}^d}$ fulfills the assumptions of Theorem 1.12 of \cite[p.\ 178]{BS}.  The collection
\[Z_n(j):= \frac{1}{n^d} \sum_{k_1,\dots, k_d=1}^n X(j_1+\tfrac{k_1}{n},\dots, j_d+\tfrac{k_d}{n})- \mathbb{E}\, X(0), \, j\in \mathbb{Z}^d,\]
is $BL(\theta')$-dependent for any $n\in\mathbb{N}$, where $\theta'_r:=\theta_{r-d}$. Indeed, let $I,J\subseteq\mathbb{Z}^d$ and let $f:\R{\# I}\to\R{}$ and $g:\R{\# J}\to\R{}$ be bounded Lipschitz functions. Put $\tilde I = I+\{1/n,2/n,\dots, 1\}^d$, $\tilde J = J+\{1/n,2/n,\dots, 1\}^d$, 
\[\tilde f: \R{\#I\cdot n^d}\to\R{}, \, (x_{1,1},\dots, x_{\#I,n^d}) \to f\Big(\frac{1}{n^d} \sum_{\ell=1}^{n^d} x_{1,\ell} - \mathbb{E} X(0), \dots, \frac{1}{n^d} \sum_{\ell=1}^{n^d} x_{\#I,\ell} - \mathbb{E} X(0)\Big)\]
and 
\[\tilde g: \R{\#J\cdot n^d}\to\R{}, \, (x_{1,1},\dots, x_{\#J,n^d}) \to g\Big(\frac{1}{n^d} \sum_{\ell=1}^{n^d} x_{1,\ell} - \mathbb{E} X(0), \dots, \frac{1}{n^d} \sum_{\ell=1}^{n^d} x_{\#J,\ell} - \mathbb{E} X(0)\Big).\]
Then we have $f(Z_{n, I}) = \tilde f(X_{\tilde I})$, $g(Z_{n, J}) = \tilde g(X_{\tilde J})$, $\Lip(\tilde f)=\Lip(f)/n^d$, $\Lip(\tilde g)=\Lip(g)/n^d$ and $\dist(\tilde I, \tilde J) \ge \dist(I,J) - d$. So
\begin{align*}
\Cov(f(Z_{n, I}),g(Z_{n, J})) & = \Cov( \tilde f(X_{\tilde I}),\tilde g(X_{\tilde J}))\\
&\le \min\{\#I\cdot n^d, \#J\cdot n^d\} \Lip(\tilde f) \Lip(\tilde g)n^d\theta_{r-d}\\
&= \min\{\#I, \#J\} \Lip(f) \Lip(g)\theta'_{r}.
\end{align*}

By Lemma \ref{L:BLlim}, the field $(Z(j))_{j\in \mathbb{Z}^d}$ is $BL(\theta')$-dependent if we can show that the finite-dimensional distributions of $(Z_n(j))_{j\in\mathbb{Z}^d}$ converge to those of $(Z(j))_{j\in \mathbb{Z}^d}$. First we will show 
\begin{equation} \lim_{n\to\infty}\mathbb{E}\, |Z_n(j)-Z(j)|=0, \quad j\in\mathbb{Z}^d.\label{eq:Zcon}\end{equation}
Let $\epsilon>0$. Due to Lemma \ref{L:contin}, the field $(X(t))_{t\in\mathbb{R}^d}$ is continuous in $1$-mean and hence there is $n$ such that 
\[\mathbb{E}\, |X(0)-X(t)| <\epsilon \mbox{ for all } t\in[0,\tfrac{1}{n}]^d. \]
Since $(X_t)_{t\in \R{d}}$ is stationary, this implies
\[\mathbb{E}\, |X(j_1+\tfrac{k_1}{n}, \dots, j_d+\tfrac{k_d}{n})-X(t)| <\epsilon \mbox{ for all } t\in[j_1+\tfrac{k_1-1}{n},j_1+\tfrac{k_1}{n}] \times \dots \times [j_d+\tfrac{k_d-1}{n},j_d+\tfrac{k_d}{n}].\]
Hence $\mathbb{E}\, |Z_n(j)-Z(j)|<\epsilon$, which finishes the proof of (\ref{eq:Zcon}).

Now let $j^{(1)}, \dots, j^{(r)}\in \mathbb{Z}^d$ and let $\delta>0$. From Markov's inequality we get 
\[\mathbb{P}\Big(\sum_{l=1}^r |Z_n(j^{(l)}) - Z(j^{(l)})|>\delta\Big) \le \frac{\sum_{l=1}^r\mathbb{E}\, |Z_n(j^{(l)}) - Z(j^{(l)})| }{\delta} \to 0, \quad n\to \infty.\]
So the finite-dimensional distributions of $(Z_n(j))_{j\in\mathbb{Z}^d}$ converge to those of $(Z(j))_{j\in \mathbb{Z}^d}$ and hence $(Z(j))_{j\in\mathbb{Z}^d}$ is $BL(\theta)$-dependent. 

By Lemma \ref{L:sameint}, the assumption that $(X(t))_{t\in \R{d}}$ is stationary implies that $(Z(j))_{j\in\mathbb{Z}^d}$ is stationary. Moreover, $(Z(j))_{j\in\mathbb{Z}^d}$ is centered, since
\[\mathbb{E}\, Z(0)=\mathbb{E}\,\int_{[0,1)^d} X(t)\, dt-\mathbb{E}\, X(0)=\int_{[0,1)^d}\mathbb{E}\, X(t)\, dt-\mathbb{E}\, X(0)=0. \]

Further,
\begin{align*}
\sum_{j\in\mathbb{Z}^d} \Cov\big(Z(0), Z(j)\big) &= \sum_{j\in\mathbb{Z}^d}\int_{[0,1)^d} \int_{j+[0,1)^d} \Cov\big(X(s), X(t)\big)\, dt\, ds\\
 &= \int_{[0,1)^d} \int_{\mathbb{R}^d} \Cov\big(X(0), X(t-s)\big)\, dt\, ds\\
 &= \int_{[0,1)^d} \int_{\mathbb{R}^d} \Cov\big(X(0), X(t)\big)\, dt\, ds\\
 &=  \int_{\mathbb{R}^d} \Cov\big(X(0), X(t)\big)\, dt.
\end{align*}

We put $Q_n:=\{j\in \mathbb{Z}^d \mid j + [0,1)^d\subseteq W_n\}$ and  $W_n^-:=\bigcup_{j\in Q_n} \big( j+[0,1)^d \big)$. As explained in the proof of \cite[Theorem 1.2]{BST}, the assumption that $(W_n)_{n\in \mathbb{N}}$ is VH-growing implies that $(Q_n)_{n\in\mathbb{N}}$ is regular growing. Now Theorem 1.12 of \cite[p.\ 178]{BS} implies that
\[\frac{\int_{W_n^-} X(t)\, dt - \lambda_d(W_n^-) \mathbb{E}\, X(0)}{\sqrt{\lambda_d(W_n^-)}} = \frac{\sum_{j\in Q_n} Z(j)}{\sqrt{\# Q_n}}\to \mathcal{N}(0,\sigma^2), \quad n\to\infty.\]

If we can show that 
\begin{equation}\frac{\int_{W_n\setminus W_n^-} X(t)\, dt - \lambda_d(W_n\setminus W_n^-)\mathbb{E}\, X(0)}{\sqrt{\lambda_d(W_n)}} \stackrel{P}{\to} 0, \quad n\to\infty, \label{eq:bd1}\end{equation}
then Slutzki's theorem will imply the assertion, since, clearly, $\sqrt{\lambda_d(W_n^-)}/\sqrt{\lambda_d(W_n)}\to 1$. 
%We put $\tilde Q_n:=\{j\in \mathbb{Z}^d \mid z + [0,1)^d\cap W_n\ne \emptyset\}$ and $W_n^+:=\bigcup_{j\in \tilde Q_n} j+[0,1)^d$. Then
We get
\begin{align*} 
\Var\Big(\int_{W_n\setminus W_n^-} X(t)\,dt\Big) &= \int_{W_n\setminus W_n^-} \int_{W_n\setminus W_n^-} \Cov(X(s), X(t))\, dt \, ds \\
&\le \int_{W_n\setminus W_n^-} \int_{\R{d}} |\Cov(X(s), X(t))|\, dt \, ds \\
&= \lambda_d(W_n \setminus W_n^-) \int_{\R{d}} |\Cov(X(0), X(t))|\, dt.
\end{align*}

Since $(W_n)_{n\in\mathbb{N}}$ is VH-growing, we get
\[\Var\bigg(\frac{\int_{W_n\setminus W_n^-} X(t)\, dt}{\sqrt{\lambda_d(W_n)}}\bigg) = \frac{\Var\big(\int_{W_n\setminus W_n^-} X(t)\, dt\big)}{\lambda_d(W_n)} \to 0, \quad n\to\infty. \]
 By the Chebyshev inequality this implies (\ref{eq:bd1}). \qed

\section{The multivariate CLT}\label{s:Multi}

In this section we extend Theorem \ref{T:uCLT} in various ways to multivariate central limit theorems.

\begin{theorem}\label{T:multCLT}
Let $\theta=(\theta_r)_{r\in\mathbb{N}}$ be a monotonically decreasing zero sequence. Let $(X(t))_{t\in\mathbb{R}^d}$ be an $\mathbb{R}^s$-valued random field. Assume that $(X(t))_{t\in\mathbb{R}^d}$ is stationary, measurable, BL($\theta$)-dependent and fulfills
\[ \int_{\mathbb{R}^d} |\Cov(X_i(0), X_j(t))| \, dt < \infty, \quad i,j=1, \dots, s.\]
Let $(W_n)_{n\in \mathbb{N}}$ be a VH-growing sequence of subsets of $\R{d}$. 
Then
\[\Big(\frac{ \int_{W_n} X_1(t)\, dt  - \mathbb{E}\, X_1(0) \lambda_d(W_n)}{\sqrt{\lambda_d(W_n)}}, \dots, \frac{ \int_{W_n} X_s(t)\, dt  - \mathbb{E}\, X_s(0) \lambda_d(W_n)}{\sqrt{\lambda_d(W_n)}} \Big)\to \mathcal{N}(0,\Sigma), \quad n\to\infty,\]
in distribution, where $\Sigma$ is the matrix with entries
\[ \int_{\mathbb{R}^d} \Cov(X_i(0), X_j(t))\, dt , \quad i,j=1, \dots, s.\]
\end{theorem}

\prf Let $u=(u_1,\dots, u_s) \in \mathbb{R}^s$. Then $(\langle X(t), u \rangle)_{t\in \R{d}}$ is $BL(\theta')$-dependent for a monotonically decreasing sequence $\theta'=(\theta'_r)_{r\in \mathbb{N}}$ with $\lim_{r \to \infty} \theta'_r =0$ due to Lemma \ref{L:LipBL}. Obviously, $(\langle X(t), u \rangle)_{t\in \R{d}}$ is stationary and measurable. We have
\[ \int_{\mathbb{R}^d} \Cov(\langle X(0), u \rangle, \langle X(t), u \rangle)\, dt =\sum_{i=1}^s\sum_{j=1}^s u_i u_j \int_{\mathbb{R}^d} \Cov( X_i(0), X_j(t))\, dt = u^T\Sigma u.\]
In particular, the integral is defined. So Theorem \ref{T:uCLT} implies
\begin{align*}
\Big\langle \Big(&\frac{ \int_{W_n} X_1(t)\, dt  - \mathbb{E}\, X_1(0) \lambda_d(W_n)}{\sqrt{\lambda_d(W_n)}}, \dots, \frac{ \int_{W_n} X_s(t)\, dt  - \mathbb{E}\, X_s(0) \lambda_d(W_n)}{\sqrt{\lambda_d(W_n)}} \Big), u \Big \rangle \\
&=\frac{ \int_{W_n} \langle X(t), u \rangle \, dt  - \mathbb{E}\, \langle X(0), u \rangle  \lambda_d(W_n)}{\sqrt{\lambda_d(W_n)}}
\to \mathcal{N}(0, u^T\Sigma u),\quad n\to\infty.
\end{align*}
Since $\langle Y, u \rangle \sim \mathcal{N}(0, u^T\Sigma u)$ for a random vector $Y\sim \mathcal{N}(0, \Sigma)$, the Theorem of Cram\'er and Wold implies the assertion. \qed\medskip

\begin{kor}\label{K:fiX}
Let $(X(t))_{t\in\mathbb{R}^d}$ be a stationary, measurable $\mathbb{R}$-valued random field and let $f_1,\dots, f_s:\mathbb{R} \to \mathbb{R}$ be functions.  %Let  Lipschitz continuous maps such that
%\[ \int_{\mathbb{R}^d} \big|\Cov\big(f_i(X(0)), f_j(X(t))\big)\big|\, dt < \infty, \quad i,j=1, \dots, s.\]
Let $(W_n)_{n\in \mathbb{N}}$ be a VH-growing sequence of subsets of $\R{d}$. 
Assume that one of the following conditions holds:
\begin{enumerate}[(i)]
\item The field $(X(t))_{t\in\mathbb{R}^d}$ is $BL(\theta)$-dependent for a monotonically decreasing zero sequence $\theta=(\theta_r)_{r\in \mathbb{N}}$, the maps $f_1,\dots, f_s$  are Lipschitz continuous and 
\[ \int_{\mathbb{R}^d} \big|\Cov\big(f_i(X(0)), f_j(X(t))\big)\big|\, dt < \infty, \quad i,j=1, \dots, s.\]
%\item The field $(X(t))_{t\in\mathbb{R}^d}$ is QA with $\mathbb{E}\, X(0)^2<\infty$. Assume that for each $i=1,\dots,s$, there is a sequence $(f_{i,n})_{n\in\mathbb{N}}$ of Lipschitz continuous functions converging pointwise to $f$ such that $(|f_{i,n}(x)|)_{n\in\mathbb{N}}$ is monotonically increasing for each fixed $x\in\mathbb{R}$. Assume, moreover, that there are $c>0$ and $\epsilon>0$ with
%\begin{equation} \Cov( X(0), X(t) ) \le c \cdot \|t\|^{-d-\epsilon}_\infty\label{e:Cov_order} \end{equation}
%for $t\in \R{d}$. 
\item The field $(X(t))_{t\in\mathbb{R}^d}$ is QA and there are $c>0$ and $\epsilon>0$ with
\begin{equation} \Cov( X(0), X(t) ) \le c \cdot \|t\|^{-d-\epsilon}_\infty, \quad t\in\mathbb{R}^d.\label{e:Cov_order} \end{equation}
The maps $f_1, \dots, f_s$ are Lipschitz continuous. 
\item The field $(X(t))_{t\in\mathbb{R}^d}$ is PA with $\mathbb{E}\, X(0)^2<\infty$. The maps $f_1,\dots, f_s$ are of locally bounded variation with $\mathbb{E}[h_{f_i}(X(0))^2]<\infty$, $i=1,\dots,s,$ and 
 there are $c>0$ and $\epsilon>0$ with
\begin{equation} \Cov\big( h_{f_i}(X(0)), h_{f_j}(X(t)) \big) \le c \cdot \|t\|^{-d-\epsilon}_\infty, \qquad  t\in \R{d}, \ i,j=1,\dots, s.\label{e:Cov_order3} \end{equation}
\end{enumerate}
Then
\[\Big(\frac{ \int_{W_n} f_1(X(t))\, dt  - \mathbb{E}\, f_1(X(0)) \lambda_d(W_n)}{\sqrt{\lambda_d(W_n)}}, \dots, \frac{ \int_{W_n} f_s(X(t))\, dt  - \mathbb{E}\, f_s(X(0)) \lambda_d(W_n)}{\sqrt{\lambda_d(W_n)}} \Big)\to \mathcal{N}(0,\Sigma), \]
as $n\to\infty$ in distribution, where $\Sigma$ is the matrix with entries
\[ \int_{\mathbb{R}^d} \Cov\big(f_i(X(0)), f_j(X(t))\big)\, dt , \quad i,j=1, \dots, s.\]
\end{kor}

Part (i) of this corollary is an immediate consequence of Lemma \ref{L:LipBL} and Theorem \ref{T:multCLT}.

\prf[ of Corollary \ref{K:fiX}(ii)]
%First we proof the assertion in the special case that the functions $f_1\dots, f_s$ are Lipschitz continuous themselves. 
The field $(X(t))_{t\in\mathbb{R}^d}$ is $BL(\theta)$-dependent by Lemma \ref{l:QABl} and thus Lemma \ref{L:LipBL} implies that the field $(f_1(X(t)), \dots, f_s(X(t)))_{t\in\R{d}}$ is also $BL(\theta)$-dependent. 

 In order to check the integrability assumptions from part (i), we put 
\[f_j^{(N)}:x\mapsto \begin{cases} -N & \mbox{if } f_j(x) < -N\\
																		f_j(x) & \mbox{if } f_j(x)\in [-N,N]\\
																		N & \mbox{if } f_j(x) > N. \end{cases}\]
Since $(X(t))_{t\in\R{d}}$ is QA, we get
\begin{align*}
| \Cov\big(f_i^{(N)}(X(0)), f_j^{(N)}(X(t))\big)| & \le \Lip(f_i^{(N)}) \cdot \Lip(f_j^{(N)}) \cdot |\Cov(X(0), X(t))| \\
& \le   \Lip(f_i) \cdot \Lip(f_j) \cdot |\Cov(X(0), X(t))|. 
\end{align*}
By the Monotone Convergence Theorem, applied to both summands of $\mathbb{E}[f_i^{(N)}(X(0))f_j^{(N)}(X(t))] - \mathbb{E}[f_i^{(N)}(X(0))] \cdot \mathbb{E}[f_j^{(N)}(X(t))]$, this yields
 \[ | \Cov\big(f_i(X(0)), f_j(X(t))\big)| \le   \Lip(f_i) \cdot \Lip(f_j) \cdot |\Cov(X(0), X(t))| . \]
%in view of \eqref{Cov_order} for an appropriate constant $\tilde c$. Since t is QA by cite {\bf REferenz}, Lemma  implies that it  .  
Moreover,  (\ref{e:Cov_order}) implies 
\[ \int_{\mathbb{R}^d} \big|\Cov\big(X(0), X(t)\big)\big|\, dt < \infty\] 
and hence
\[ \int_{\mathbb{R}^d} \big|\Cov\big(f_i(X(0)), f_j(X(t))\big)\big|\, dt < \infty, \quad i,j=1, \dots, s. \]
So part (i) yields the assertion. %, in the special case that $f_1,\dots f_s$ are Lipschitz continuous.

 \qed\medskip

\prf[ of Corollary \ref{K:fiX}(iii)] 
Since $(X(t))_{t\in\mathbb{R}^d}$ is PA, the random field $(h_{f_1}(X(t)), \dots, h_{f_s}(X(t)))_{t\in\mathbb{R}^d}$ is also PA, see Theorem 1.8(d) of \cite[p.\ 7]{BS}, and therefore QA. By Lemma \ref{l:QABl} it is $BL(\theta)$-dependent for some monotonically decreasing zero sequence $\theta$. Hence $(f_1(X(t)), \dots, f_s(X(t)))_{t\in\mathbb{R}^d}$ is $BL(\theta')$-dependent for some monotonically decreasing zero sequence $\theta'$ by Lemma \ref{L:bv_decomp} and Lemma \ref{L:LipBL}. 

Clearly, the field $(f_1(X(t)), \dots, f_s(X(t)))_{t\in\mathbb{R}^d}$ is also stationary and measurable.  

%Since $(f_1(X(t)), \dots, f_s(X(t)))_{t\in\mathbb{R}^d}$ is measurable and stationary, it is stochastically continuous, see \cite[Prop.\ 3.1]{Roy}. 

Moreover, \eqref{e:Cov_order3} implies
\[ \int_{\mathbb{R}^d} \big|\Cov\big(h_{f_i}(X(0)), h_{f_j}(X(t))\big)\big|\, dt < \infty, \quad i,j=1, \dots, s. \]
Now Lemma \ref{L:bv_decomp} and the QA property of $(h_{f_1}(X(t)), \dots, h_{f_s}(X(t)))_{t\in\mathbb{R}^d}$ give
\[ \int_{\mathbb{R}^d} \big|\Cov\big(f_i(X(0)), f_j(X(t))\big)\big|\, dt \le \int_{\mathbb{R}^d} \big|\Cov\big(h_{f_i}(X(0)), h_{f_j}(X(t))\big)\big|\, dt < \infty, \quad i,j=1, \dots, s. \]
%can be obtained exactly the same way as in the proof of part (ii).

So Theorem \ref{T:multCLT} yields the assertion.\qed

\end{document}